\documentclass[final]{amsart}

\usepackage{color}
\usepackage{newlfont}

\usepackage[english, status=draft]{fixme}
\fxusetheme{color}
\usepackage{comment}
\usepackage{pictexwd,dcpic}
\usepackage{amssymb}
\usepackage{amsmath, amscd}
\usepackage{latexsym}
\usepackage{amsthm}
\usepackage{slashed}
\usepackage{mathtools}
\usepackage[toc,page]{appendix}
\usepackage[utf8]{inputenc}
\usepackage[T1]{fontenc}
\usepackage{fixme}
\usepackage{slashed}

\makeatletter
\DeclareFontFamily{OMX}{MnSymbolE}{}
\DeclareSymbolFont{MnLargeSymbols}{OMX}{MnSymbolE}{m}{n}
\SetSymbolFont{MnLargeSymbols}{bold}{OMX}{MnSymbolE}{b}{n}
\DeclareFontShape{OMX}{MnSymbolE}{m}{n}{
	<-6>  MnSymbolE5
	<6-7>  MnSymbolE6
	<7-8>  MnSymbolE7
	<8-9>  MnSymbolE8
	<9-10> MnSymbolE9
	<10-12> MnSymbolE10
	<12->   MnSymbolE12
}{}
\DeclareFontShape{OMX}{MnSymbolE}{b}{n}{
	<-6>  MnSymbolE-Bold5
	<6-7>  MnSymbolE-Bold6
	<7-8>  MnSymbolE-Bold7
	<8-9>  MnSymbolE-Bold8
	<9-10> MnSymbolE-Bold9
	<10-12> MnSymbolE-Bold10
	<12->   MnSymbolE-Bold12
}{}

\let\llangle\@undefined
\let\rrangle\@undefined
\DeclareMathDelimiter{\llangle}{\mathopen}%
{MnLargeSymbols}{'164}{MnLargeSymbols}{'164}
\DeclareMathDelimiter{\rrangle}{\mathclose}%
{MnLargeSymbols}{'171}{MnLargeSymbols}{'171}
\makeatother

\newcommand{\boldnabla}{\mbox{\boldmath$\nabla$}}

\newcommand{\SqInt}{\operatorname{L}^2}
\newcommand{\diff}{\mathrm{d}}

\newcommand{\Half}{\mathbb H}

\newcommand{\conj}[1]{\overline{#1}}

\newcommand{\Real}{\mathbb{R}}

\newcommand{\Integer}{\mathbb{Z}}

\newcommand{\Complex}{\mathbb{C}}

\newcommand{\Torus}{\mathbb{T}}

\newcommand{\Hilb}{ {\mathcal H}}

\newcommand{\Line}{\mathcal{L}}
\newcommand{\Cont}{\operatorname{C}}
\newcommand{\Smooth}{\Cont^{\infty}}

\newcommand{\bpar}{\mathrm{b}}

\newcommand{\rprod}[2]{\left(#1\text{,}#2\right)}
\newcommand{\krprod}[2]{\left(#1\text{,}#2\right)_k}
\newcommand{\fker}[2]{\left\langle#1\text{,}#2\right\rangle}
\newcommand{\inn}[2]{\llangle#1\text{,}#2\rrangle}

\newcommand{\comm}[1]{\left[ #1\right]}

\newcommand{\opG}{\mathsf G}

\newcommand{\opF}{\mathsf{F}}

\newcommand{\opb}{\mathsf b}

\newcommand{\opf}{\mathsf f}
\newcommand{\opg}{\mathsf g}
\newcommand{\Zcal}{\mathcal Z}
\newcommand{\kZcal}{\Zcal_{(k)}}

\newcommand{\p}{\prime}
\newcommand{\pp}{{\prime\prime}}

\newcommand{\glie}{\mathfrak{g}}

\newcommand{\torus}{\mathfrak{t}}

\newcommand{\coroot}{\Lambda^R}
\newcommand{\weights}{\Lambda^w}
\newcommand{\kweights}{\Lambda^w_{(k)}}

\theoremstyle{definition}

\newtheorem{remark}{Remark}

\newtheorem{proposition}{Proposition}
\newtheorem{lemma}{Lemma}
\newtheorem{theorem}{Theorem}

\DeclareMathOperator{\SU}{SU} 
\newcommand{\im}{\mathop{\fam0 Im}\nolimits}
\newcommand{\spane}{\mathop{\fam0 Span}\nolimits}
\newcommand{\tr}{\mathop{\fam0 Tr}\nolimits}
\newcommand{\re}{\mathop{\fam0 Re}\nolimits}

\newcommand{\Hom}{\mathop{\fam0 Hom}\nolimits}

\newcommand{\vol}{\mathop{\fam0 Vol}\nolimits}

\newcommand{\id}{\mathop{\fam0 Id}\nolimits}

\newcommand{\ra}{\mathop{\fam0 \rightarrow}\nolimits}

\newcommand{\cH}{ {\mathcal H}}

\renewcommand{\epsilon}{\varepsilon}

\renewcommand{\phi}{\varphi}

\newcommand{\abs}[1]{\lvert #1 \rvert}

\newcommand{\setZ}{\mathbb{Z}}
\newcommand{\setR}{\mathbb{R}} 
\newcommand{\setC}{\mathbb{C}}

\newcommand{\Vol}{\mathrm{Vol}}

\newcommand{\dvol}{\diff\text{\em vol}}

\renewcommand{\Re}{\mathrm{Re}}

\def\XXint#1#2#3{{\setbox0=\hbox{$#1{#2#3}{\int}$}
\vcenter{\hbox{$#2#3$}}\kern-.5\wd0}}

\newcommand{\mT}{\mathbb{T}}

\newcommand{\rhoq}{\eta}

\begin{document}

\title[Genus one complex Quantum Chern-Simons represenation]{The genus one Complex Quantum Chern-Simons representation of the Mapping Class Group}
\author{J{\o}rgen Ellegaard Andersen}
\address{Center for Quantum Geometry of Moduli Spaces\\
        University of Aarhus\\
        DK-8000, Denmark}
\email{andersen@qgm.au.dk}

\author{Simone Marzioni}
\address{Center for Quantum Geometry of Moduli Spaces\\
	University of Aarhus\\
	DK-8000, Denmark}
\email{marzioni@qgm.au.dk}

\thanks{Supported in part by the center of excellence grant ``Center for quantum geometry of Moduli Spaces" from the Danish National Research Foundation (DNRF95)}

\begin{abstract}
In this paper we compute explicitly, following Witten's prescription, the quantum representation of the mapping class group in genus one for complex quantum Chern-Simons theory associated to any simple and simply connected complex gauge group $G_\Complex$. We use a generalization of the Weil-Gel'fand-Zak transform to exhibit an explicit expression for the representation. 
\end{abstract}

\maketitle

\section{Introduction}
In his seminal paper \cite{W1}, Witten introduce the study of quantum Chern-Simons theory with non-abelian gauge group $G$. For the compact gauge groups $G=SU(N)$, the so called Witten-Reshetikhin-Turaev TQFT's was first constructed by Reshetikhin and Turaev for $G=SU(N)$ \cite{RT1,RT2,T} and subsequently by Blanchet, Habegger, Masbaum and Vogel in \cite{BHMV1,BHMV2,B1} using skein theory. Recently Ueno jointly with the first author of this paper showed in a series of four papers \cite{AU1,AU2,AU3,AU4} that the Witten-Reshetikhin-Turaev TQFT is the same as the TQFT coming from Conformal Field Theory \cite{TUY, BK} as was also proposed to be the case in Witten's original paper \cite{W1}. Witten further suggested in the same paper that the geometric quantization of the moduli space of flat connections should be related to this theory and he developed this approach further in his joint paper with Axelrod and Della Pietra \cite{ADW}. Following shortly after, Hitchin gave a rigorous account of this work in \cite{H}. For a purely differential geometric account of the construction of this connection see \cite{A9, AG1, AGL}. By combining the work of Laszlo \cite{La1} with the above mention work of Ueno and the first author, it has now been confirmed that one can use the geometric quantization of the moduli space of flat connections as an alternative construction of the Witten-Reshetikhin-Turaev TQFT. This has been exploited in, among others, the works \cite{A6, AMU, A7, AHi, A8, A10,  1408.2499}.

In the paper \cite{W2} Witten also proposed a way to construct the mapping class groups representations of quantum Chern-Simons theory for complex gauge groups. This theory received much less attention, but the time now seems ripe to fully develop this theory. Recently Kashaev has, jointly with the first author of this paper, rigorously constructed the $n=2$ theory in \cite{AK1,AK2, AK3}. First the level $k=1$ theory was examined in detail in \cite{AK1,AK2} and recently has the picture for all levels become clear and in fact, a very general scheme has been proposed in \cite{AK3} for construction of TQFT's from Pontryagin self-dual Locally compact groups. It is further verified in that paper that the scheme for the LCA $\setR\times \setZ_k$ yields quantum Chern-Simons theory for $PSL(2,\setC)$ at level $k$. See also \cite{AS} for further discussion of the general level theory. This should be seen in parallel to the developments on the index by Garoufalidis and Dimofte  in \cite{1208.1663,MR3073925}, which should be related to the level $k=0$ theory and further their work on the quantum modularity conjecture \cite{DGa2}. In the physics literature, the complex quantum Chern-Simons theory has been discussed from a path integral point of view in a number of
papers \cite{MR3250765,BNW, MR3148093,MR2551896,MR3080552,MR2465747,MR2134725,Hik1,Hik2,W3} (see also references in these) and latest by
Dimofte \cite{Dimofte3d3d, Di3} using the more advanced 3d-3d correspondence.

In parallel to this Gammelgaard in collaboration with the first author of this paper, has given \cite{AG2} a finite dimensional differential geometric proof of the projective flatness of the connections proposed by Witten in \cite{W2}. The basic setup is here to consider the space of all sections of the Chern-Simons line bundle over the moduli space of flat $G$ connections as the polarised Hilbert space for the quantisation of the $G_\Complex$ moduli space, where the polarisation on the $G_\Complex$ moduli space depends on the choice of a complex structure on the underlying surface, just as in the case of the complex polarizations of the moduli space for the compact group $G$. Witten shows that one also in this non-compact case gets a connection, which he argues using partially physics arguments, should be flat. In the work \cite{AG2} the first author and Gammelgaard gave a differential geometric construction of this connection, which they named the {\em Hitchin-Witten} connection. Further they gave a proof that this Hitchin-Witten connection is projectively flat for higher genus surfaces. 

For the genus one case, it is immediately clear that Witten's proposal produces a flat connection. The purpose of this paper is to compute the resulting representation of the mapping class group in this case. Let us now review the setup from Witten's paper\footnote{At this point, we wish however, to remark that when considering the moduli space, the non-stable bundles are ignored. As can be seen from \cite{AGP}, working with the stack of all bundles as opposed to the moduli spaces of semi-stable bundles, might effect significantly the answer exactly in genus one.}.

Suppose $G$ is a simple, simply connected compact Lie group and let $G_{\Complex}$ be its complexification. Let $\Sigma$ be a closed oriented surface of genus one. Then we get the following description of the moduli space $M_G$ of flat $G$ connection on  $\Sigma$
$$
  M_G \cong \mT_G\times \mT_G/W,
$$ 
where $\mT_G$ is the maxim torus of $G$ and $W$ is the Weyl group of $G$. We thus have the projection map
$$ \pi : \mT_G\times \mT_G\to M_G,$$
whose fiber over a generic point in the target is a copy of $W$, but since the action of $W$ on $\mT_G\times \mT_G$ has fixed points, this is not a global covering map.
Using the Chern-Simons action functional and the correct normalisation of the invariant inner product on the Lie algebra of $G$ (see e.g. \cite{Fr}, \cite{RSW}), one can construct a complex line bundle $ \Line_G$ over $ M_G $ and from the construction we get the following identification for all integers $k$
\begin{align}
  \nonumber
  C^\infty(M_G, \Line_G^k) \cong C^\infty(\mT_G\times \mT_G, \pi^* \Line_G^k)^W.
\end{align}
Introducing $\cH^{(k)} \equiv C^\infty(M_G, \Line_G^k)$ and  $\widetilde\cH^{(k)}\equiv C^\infty(\mT_G\times \mT_G, \pi^*\Line_G^{k})$, we of course have that
$$ \cH^{(k)} = (\widetilde\cH^{(k)})^W.$$

Consider ${\Half}$, the upper half plane, which is the  Teichm\"{u}ller space for $\Sigma$, in the sense that it parametrises all marked complex structures on $\Sigma$ up to natural equivalance (see also \cite{A5} for the general torus case). 
The construction of Complex Quantum Chern-Simons Theory with gauge group $G_{\Complex}$ outlined in \cite{W2} provides an identification of the Quantum Hilbert space of this theory with the pre-Quantum Hilbert space for the theory relative to the compact group $G$, that is $\cH^{(k)}$. However, such an identification depends on a choice of complex structure on the $\Sigma$. The complex structure on $\Sigma$ induces a real polarization on the moduli space $M_{G^{\mathbb C}}$, the moduli space of flat $G^{\mathbb C}$-connection, via the Hodge star operator on $\Sigma$ as described in \cite{W2}. Consequently we have a trivial vector bundle $\cH^{(k)}\times\Half\rightarrow \Half$ with a non-trivial flat Hitchin-Witten connection.
For the genus one moduli spaces the Hitchin--Witten connection $\boldnabla$ happens to have a simple description as follows. Fix a complex number $t= k+is$, such that $k$, the real part of $t$, is an integer. Let $\diff$ be the trivial connection on the bundle $\Hilb^{(k)}\times \Half$. As Witten argues in \cite{W2}, if one chooses $r$ such that
\begin{align}
\nonumber
e^{-4r} &= -\frac{k-is}{k+is},
\end{align}
then one finds that
\begin{align}\label{eq:HWconj}
e^{-r\Delta}\circ\diff\circ e^{r\Delta} &= \boldnabla,
\end{align}
where $\Delta = \Delta_\sigma$ is the Laplace operator on $M_G$, which of course depends on $\sigma\in\Half$ parametrising the family of complex structures compatible with the Goldman symplectic form on $M_G$, which are induced from the natural family of complex structures on $\Sigma$.
This means that parallel sections of $ \mathcal{H}^{(k)} \times \Half$
are of the form
\begin{align}
\nonumber
s_\sigma &= e^{-r\Delta_\sigma} s_0
\end{align}
for any $s_0 \in C^\infty(M_G, \Line^k)$.

The mapping class group $\Gamma = SL(2, \setZ)$  acts on $M_G$ and classical Chern-Simons theory allows one to lift this action to $\Line^k$, thus one gets an induced pre-quantum action 
$$\rho_k : \Gamma \ra B(\cH^{(k)}),$$
where $B(\cH^{(k)})$ refers to the continuous operators which respect to the Frechet topology on $\cH^{(k)}$.
We call this action the \textit{pre-Quantum} \textit{representation} of $\Gamma$. Note that we can also complete this to an action on the Hilbert space of $L_2$-sections of $L^k$ over $M_G$, in which case $\Gamma$ acts by bounded operators on this Hilbert space.

What we are seeking, however, are the \textit{quantum} \textit{representations} $\rhoq_{t}$ which arise when we take into account the full action of $\Gamma$ on the bundle $\cH^{(k)}\times \Half$ composed with Hitchin-Witten parallel transport $\mathcal{P}_{\sigma_0,\sigma_1} \colon  \left(\cH^{(k)}, \sigma_1\right) \ra \left(\cH^{(k)}, \sigma_0\right)$. That  is, we first consider the action
\begin{align}
\tilde{\phi}\colon \left(\cH^{(k)}, \sigma\right) \ra \left(\cH^{(k)}, \phi_*\sigma\right)
\end{align}
and then compose it with the parallel transport, to get a representation on $\left(\cH^{(k)}, \sigma\right)$
\begin{align}
\rhoq_{t,\sigma}(\phi) \equiv \mathcal{P}_{\sigma,\phi^*\sigma} \circ \tilde \phi :\left(\cH^{(k)}, \sigma\right) \ra \left(\cH^{(k)}, \sigma\right).
\end{align}
However, the description \eqref{eq:HWconj} of $\boldnabla$ leads to the following identification
\begin{proposition}\label{pr:qr}
	The densely defined transformation 
	$$ e^{r\Delta_\sigma} : \cH^{(k)} \ra \cH^{(k)}$$
	conjugates $\rho_k$ to $\rhoq_{t,\sigma}$ for all $\sigma \in \Half$. Thus for all $\phi\in\Gamma$, we have that
	\begin{align}
	\rhoq_{t, \sigma}(\phi) = e^{-r\Delta_{\sigma}}\circ\rho_{k}(\phi)\circ e^{r\Delta_\sigma}
	\end{align}
	on the domain of $e^{r\Delta_\sigma}$.
\end{proposition}
Let us consider the following generators of $\Gamma$
\begin{align}
  \nonumber
  S &=
  \begin{pmatrix}
    0 & -1\\ 1 & 0
  \end{pmatrix}
  &
  T &=
  \begin{pmatrix}
    1 & 1 \\ 0 & 1
  \end{pmatrix},
\end{align}
which satisfies the relations
\begin{align}
  \nonumber
  S^2 &= (ST)^3 & S^4 &= \id.
\end{align}

In this paper we compute the representations $\rho_k$ and $\rhoq_t$ completely explicitly. We use standard in notation from Lie theory, which we briefly recall in Appendix \ref{ap:Lie}, where we further introduce a few quantities we will need in this paper. In particular we recall that $\torus\subset \glie$ is the Lie algebra of $\Torus_G$ and that we denote by $\kZcal$ the finite set of $\frac 1k$-scaled weights modulo the root lattice of $\glie$.\\
In  Section \ref{sc:WGZ} we present a generalization of the Weil-Gel'fand-Zak transform which provides an isometry
$$ Z :  \cH^{(k)} \ra \left(\SqInt(\torus)\otimes \setC^{\kZcal}\right)^{W}.$$ 
The $W$ invariance on the right is with respect to the following action of $w\in W$
\begin{align}
&w\cdot f_\gamma(\theta) = f_{w\gamma}(w\theta)  &\text{for }\gamma\in\kZcal\text{, }\theta\in\Torus.
\end{align}
For technical reasons, we further conjugate with a discrete Fourier transform  $\opF_{\kZcal}$ (see the definition from equation \eqref{eq:FiniteOps}).
We use the following notation for the corresponding pre-quantum operators action on the square integrable functions
\begin{align}\nonumber
 \hat{\rho}_k(S) &=\left(Z\circ\opF_{\Zcal}\right)^{-1} \circ\tilde S \circ \left(Z\circ\opF_{\Zcal}\right), \\
  \hat \rho_k(T) &= \left(Z\circ\opF_{\Zcal}\right)^{-1} \circ\tilde T \circ \left(Z\circ\opF_{\Zcal}\right)
\end{align}
and correspondingly for the quantum representation
\begin{align}\nonumber
  \hat{\rhoq}_t( S)& = \left(Z\circ\opF_{\Zcal}\right)^{-1} \circ\rhoq_t( S) \circ \left(Z\circ\opF_{\Zcal}\right), \\
  \hat\rhoq_t(T) &=\left(Z\circ\opF_{\Zcal}\right)^{-1} \circ\rhoq_t(T) \circ \left(Z\circ\opF_{\Zcal}\right).
\end{align}
Considering the $W$-invariant subspace  of $\SqInt(\torus)\otimes \setC^{\kZcal}$ we get the decomposition
\begin{align}
\nonumber
\Hilb^{(k)}&\simeq \left(\SqInt(\torus)\otimes\Complex^{\kZcal}\right)^W \\
\label{eq:SpaceDecomp}
&= \left(\SqInt(F_0))\otimes \overline{C}_{(k)}\right)\bigoplus\left(\SqInt(F_0)\otimes {C}_{(k)}\right)
\end{align}
for which we refer again to appendix \ref{ap:Lie} for definitions and notations. 

\begin{theorem}\label{th:pqr}
	Let $j$ and $\omega$ be complex numbers satisfying
	\begin{align}
	&j^4 =1,
	&\omega^3 = i^{n/2}j^{-1}.
	\end{align}
	The pre-Quantum representation $\hat\rho_k$ splits into the direct sum of two representations, induced by the decomposition (\ref{eq:SpaceDecomp})
	\begin{align*}
	\hat\rho_k & =  
	(\hat\rho^\p_{k,0}\otimes\hat\rho^\pp_{k,0})\oplus(\hat\rho^\p_{k,1}\otimes\hat\rho^\pp_{k,1}),
	\end{align*}
	where, for $\theta\in F_0$, $\alpha, \beta\in \overline{I}_{k}$ and $\alpha^\p, \beta^\p\in {I}_{k}$
	\begin{align*}
	\begin{split}
	&\hat\rho^\p_{k,0}(S)(f) (\theta)= j\int_{F_0}f(\tilde{\theta})\sum_{w\in W} e^{2\pi i\krprod{w\theta}{\tilde{\theta}}}\dvol_k(\tilde{\theta}),\\
	&\hat\rho^\p_{k,1}(S)(f) (\theta)= j\int_{F_0}f(\tilde{\theta}) \sum_{w\in W}\det(w) e^{2\pi i\krprod{w\theta}{\tilde{\theta}}}\dvol_k(\tilde{\theta}),\\
	&\hat\rho^\p_{k,0}(T)(f) (\theta) = \hat\rho^\p_{k,1}(T)(f) (\theta) = \omega e^{-\pi i \krprod{\theta}{\theta}}f(\theta),\\
	\end{split}\\
	&\hat\rho^\pp_{k,0}(S)_{\alpha,\beta} = \,\frac{j^{-1}}{\sqrt{\abs{\kZcal}} }
		\sum_{w\in W} \det(w)e^{2\pi i \krprod{w\alpha}{\beta}}\\
	&\hat\rho^\pp_{k,1}(S)_{\alpha^\p,\beta^\p} =\,\frac{j^{-1}}{\sqrt{\abs{\kZcal}} }
	\sum_{w\in W}e^{2\pi i \krprod{w\alpha^\p}{\beta^\p}}\\
	&\hat\rho^\pp_{k,0}(T)_{\alpha,\beta} =  \omega^{-1}e^{\pi i \krprod{\alpha}{\alpha}}\delta_{\alpha,\beta},\\
	&\hat\rho^\pp_{k,1}(S)_{\alpha^\p,\beta^\p} = \omega^{-1}e^{\pi i \krprod{\alpha^\p}{\alpha^\p}}\delta_{\alpha^\p,\beta^\p}.
	\end{align*}
\end{theorem}
Let $\bpar\in\Complex$ such that $\abs{\bpar} =1$, $\re(\bpar)>0$ and 
\begin{equation}
is = k\frac{1-\bpar^2}{1+\bpar^2}\text{,}
\end{equation}
\begin{theorem}\label{th:qr}
		The quantum representation $\hat\rhoq_t$ is obtained by conjugation with an unitary operator
		\begin{align*}
		\hat\rhoq_{t}(\phi) = e^{-r\hat\Delta_{\sigma}}\circ\hat\rho_{k}(\phi)\circ e^{r\hat\Delta_\sigma}.
		\end{align*}
		It decomposes itself as
		\begin{align*}
		\hat\rhoq_t & =  
		(\hat\rhoq^\p_{t,0}\otimes\hat\rho^\pp_{k,0})\oplus(\hat\rhoq^\p_{t,1}\otimes\hat\rho^\pp_{k,1}).
		\end{align*}
		For the particular choice of complex structure $\sigma = i\bpar$ the representations on $\SqInt(F_0)$ take the following explicit integral form
		\begin{align*}
	&\hat\rhoq^\p_{t,0}(S)(f) (\theta)= je^{\pi(\bpar-\conj{\bpar}) \krprod{\theta}{\theta}}\int_{F_0}f(\tilde{\theta})\sum_{w\in W} e^{2\pi i\krprod{w\theta}{\tilde{\theta}}}e^{-\pi(\bpar-\conj{\bpar}) \krprod{\tilde\theta}{\tilde\theta}}\dvol_k(\tilde{\theta}),\\
	&\hat\rhoq^\p_{t,1}(S)(f) (\theta)= je^{\pi(\bpar-\conj{\bpar}) \krprod{\theta}{\theta}}\int_{F_0}f(\tilde{\theta}) \sum_{w\in W}\det(w) e^{2\pi i\krprod{w\theta}{\tilde{\theta}}}e^{-\pi(\bpar-\conj{\bpar}) \krprod{\tilde\theta}{\tilde\theta}}\dvol_k(\tilde{\theta}),\\
	&\hat\rhoq^\p_{t,0}(T)(f) (\theta) =\omega i^{-\frac n2} e^{\pi(\bpar-\conj{\bpar}) \krprod{\theta}{\theta}}
	\int_{F_0} f (\tilde{\theta}) \sum_{w\in W}
		e^{{\pi i}\krprod{w\theta- \tilde{\theta}}{w\theta-\tilde{\theta}} } e^{-\pi(\bpar-\conj{\bpar}) \krprod{\tilde\theta}{\tilde\theta}}\dvol_k(\tilde \theta),\\
	& \hat\rhoq^\p_{t,1}(T)(f) (\theta) =\omega i^{-\frac n2} e^{\pi(\bpar-\conj{\bpar}) \krprod{\theta}{\theta}}
	\int_{F_0} f (\tilde{\theta})\!\! \sum_{w\in W}\det(w)
	e^{{\pi i}\krprod{w\theta- \tilde{\theta}}{w\theta-\tilde{\theta}} } e^{-\pi(\bpar-\conj{\bpar}) \krprod{\tilde\theta}{\tilde\theta}}\dvol_k(\tilde \theta),
		\end{align*}
\end{theorem}
In Section \ref{sc:HWC}, Lemma \ref{lm:Kernel}, we provide the explicit expression for $e^{-r\hat{\Delta}_\sigma}$ for generic $\sigma$.
In Section \ref{sc:Compact} we provide comparisons to the known result for quantum Chern--Simons theory with compact gauge group $G$.
	
\section{Moduli space of flat connections and the pre-quantum line bundle}
\label{sc:MS}
In this section we briefly recall the construction of the pre-quantum space for Chern-Simons theory. We refer to \cite{Fr} and \cite{ADW} for more details. See Appendix \ref{ap:Lie} for notations on Lie theory.
Given a closed surface $\Sigma$ of genus $g = 1$, its fundamental group is abelian, 
precisely $$ \pi_1(\Sigma) = \Integer\times \Integer\text{.}$$
Let $G$ be a compact, simple, connected and simply connected Lie group. The moduli space of flat connections is then
\begin{equation}
M = \Hom \left(\pi_1(\Sigma),G\right)/G = \left(\Torus\times\Torus\right)/W,
\end{equation} 
where $\Torus$ is the maximal torus of $G$ and  $W$ is the Weyl group of $G$. We can lift the description to quotient of its tangent space
\begin{align}\label{eq:Mgauge}
M = (\torus\oplus\torus)/\left((\coroot\times\coroot)\rtimes W\right).
\end{align}
Since any principal $G$-bundle on $\Sigma$ is trivial, we can fix a trivialization and so a $G$-connection $A \in \Omega(\Sigma,\glie)$ is a one form on $\Sigma$ with values in the lie algebra $\glie$. In particular, we can identify $\torus\oplus\torus$ with a subspace of flat connections with values in $\torus$. The group $(\coroot\times\coroot)\rtimes W$ corresponds to the gauge transformations preserving this subspace of flat connections inside $\Omega(\Sigma, \torus)$.
Explicitly if $x, y$ are local coordinates on $\Sigma = \Real^2/\Integer^2$, we can write 
$A = \theta_1\diff x +\theta_2 \diff y$, where $\theta_1\text{,} \, \theta_2 \in \torus$ and identify $A$ with $\theta_1\oplus\theta_2 \in\torus\oplus\torus$. 
 Recall the preferred inner product $\rprod{\cdot}{\cdot}_1$ on $\torus$ with the normalization specified in Appendix \ref{ap:Lie}.
The Atiyah-Bott \cite{AB} symplectic form  on $\torus\oplus\torus$ is given by
\begin{equation}
\omega(A,A^{\prime}) \equiv {2\pi} \int_{\Sigma} \left(A \wedge A^{\prime}\right)_1 = 2\pi\left( \rprod{\theta_1}{\theta^\p_2}_1-\rprod{\theta^\p_1}{\theta_2}_1\right)\text{,}
\end{equation}
where $A= \theta_1\oplus \theta_2$ and $A^\p = \theta^\p_1\oplus \theta_2^\p$.\\

Let $\left\{\opb_1,\dots,\opb_n\right\}$ be a basis for $\torus$.
We can write coordinates for $A= \theta_1\diff x + \theta_2\diff y$ as
\begin{align}\label{eq:thetacoor}
&\theta_1 = \sum_{j=1}^{n}u_j{\opb_j},
&\theta_2 = \sum_{j=1}^{n}v_j{\opb_j}.
\end{align}
A complex structure on $\Sigma$ is equivalent to a choice of holomorphic coordinate $w= x+\tau y$ with $\im \tau >0$. This complex structure induces a complex structure $J = - *$ on $M$, where $*$ is the hodge star operator. $J$ is compatible with $\omega$ and we can write holomorphic coordinates on $M$
\begin{align}\label{eq:holocoord}
&z_j = u_j + \sigma v_j
&\conj{z}_j = u_j + \conj{\sigma} v_j 
\end{align}  
where 
$$\sigma = - \frac 1 \tau \text.$$
The action of $w\in W$ on $\theta_1\oplus\theta_2$ is given by
\begin{align}
w(\theta_1,\theta_2) = \left(\sum_{j} u_jw({\opb_j}),  \sum_{j} v_jw({\opb_j}) \right)
\end{align}
Let $C^{(k)}= (C_{jl}^{(k)})_{j,l=1,\dots,n}$ be the matrix defined as $C_{jl}^{(k)} \equiv\krprod{{\opb_j}}{{\opb_l}}$. The symplectic form $\omega$ takes the following explicit form
\begin{align}
\omega = 2\pi \sum_{j,l} C_{jl}^{(1)} \diff u_j\wedge\diff v_l.
\end{align}

The level $1$ pre-Quantum line bundle $\Line \rightarrow M$ defined by classical Chern-Simons theory (see \cite{Je,Fr,ADW}) can be constructed starting with the trivial bundle $\widetilde{\Line}\equiv\torus\oplus\torus\times\Complex$.  We lift the action of $\coroot\times\coroot$ to $\widetilde{\Line}$ as follows, for $(A,\zeta)\in\torus\oplus\torus\times\Complex$ and $\lambda= (\lambda_1,\lambda_2)\in\coroot\times\coroot$,
\begin{align}\label{eq:LatAction}
&(A,\zeta)\cdot\lambda = (A+\lambda, e_{\lambda}(A)\zeta),\\
&\text{where } e_{\lambda}(A) \equiv (-1)^{\rprod{\lambda_1}{\lambda_2}_1}\exp-\frac{i}{2}\left(\omega(A,\lambda)\right).
\end{align} 
This defines a line bundle $\pi^*\Line\rightarrow \Torus\times\Torus$, where $\pi$ is the projection onto $M$. The first Chern class of $\pi^*\Line$ is $c_1(\pi^*\Line) = \left[\frac{\omega}{2\pi}\right]\in H^2(\Torus\times\Torus,\Integer)$.\\
To lift the action of $W$ to $\pi^*\Line$ one takes the trivial lift to $\widetilde{\Line}$
\begin{align}\label{eq:WeylLineAct}
&(A, \zeta)\cdot w = (w(A), \zeta)\text{, for every }w\in W.
\end{align}
It is simple to prove that such maps are equivariant with respect to the action of $\coroot\times\coroot$ (see \cite{Je} Lemma A.2). This completes the description of the Chern--Simons line bundle $\Line\rightarrow M$.\\
\begin{remark}
The $e_{\lambda}$ are called {\em multipliers} in the literature. The factor $(-1)^{\rprod{\lambda_1}{\lambda_2}_1}$ is usually referred as a \emph{theta--characteristic} and it is not uniquely defined by the Chern--Simons action. In fact it is irrelevant in the definition of $\Line\rightarrow M$ but a choice of it is needed to lift it to $\pi^*\Line$ (see \cite{Je}).
\end{remark}
The line bundle $\Line$ support a preferred pre-Quantum connection, given explicitly on $\tilde{\Line}$ by
\begin{align}\label{eq:pqConnection}
&\nabla = \diff + \alpha\\
&\alpha_{A}(X) = \omega(A, X)\text{, for any } X\in T_{A}(\torus\oplus\torus).
\end{align} 
This connection has curvature $F_{\nabla} = -i\omega$ and is compatible with the actions \eqref{eq:LatAction} and \eqref{eq:WeylLineAct}, therefore it descends to a well defined connection on $\Line$.\\

The line bundle $\Line^k \equiv \Line^{\otimes k}$ can be easily obtained from the description above substituting the inner product $\rprod{\cdot}{\cdot}_1$ with the $k$--scaled inner product $\krprod{\cdot}{\cdot}$.\\

Consider the space $\widetilde \Hilb^{(k)} = C^{\infty} (\Torus\times\Torus, \pi^*\Line^{k})$.
We are interested in the level $k$ quantum  space 
$$\Hilb^{(k)} \equiv C^{\infty} (M, \Line^{k}) = C^{\infty} (\Torus\times\Torus, \pi^*\Line^{k})^{W}$$
 of Weyl group invariant smooth sections of $\pi^*{\Line}^{2k}$.\\
The inner product given by
\begin{equation}
\phi\cdot\psi(p) = \phi(p)\conj{\psi(p)}\text, \qquad p\in\Torus\times\Torus\text, 
\quad \phi\text,\, \psi \in \widetilde\Hilb^{(k)},
\end{equation}
is parallel with respect to our choice of $\nabla$.
This allows us to define an inner product on $\widetilde{\Hilb}^{(k)}$ as
 \begin{equation}\label{eq:hermitianDef}
 \inn{\phi}{\psi} = \frac{1}{\abs{\Lambda}}\int_{\Torus\times\Torus}\phi\cdot\conj\psi\,\dvol_k
 \text, \qquad p\in\Torus\times\Torus\text, 
 \quad \phi\text,\, \psi \in \tilde\Hilb^{(k)}.
 \end{equation}
Here $\dvol_k$ is the volume form induced by the inner product $\rprod{\cdot}{\cdot}_k$.\\

An element $\gamma\in\Gamma$ acts by pull back on $M$ mapping the equivalence class or representations $[\rho]\in M$ to $[\rho\circ \gamma]$.
Classical Chern-Simons theory enables us to lift the action to the line bundle $\Line^k$ to get an operator 
\begin{equation}
\tilde \gamma\colon \Hilb^{(k)} \longrightarrow \Hilb^{(k)} \text,
\end{equation}
by defining, on $\tilde{\Line}$
\begin{align}
(A,\zeta)\cdot\gamma = (\gamma^*(A), \zeta).
\end{align}
It can be easily shown that this lift is equivariant with respect to $\coroot\times\coroot$ (see \cite{Je}).
This is not sufficient to define a complex quantum representation of $\Gamma$, indeed an identification of the $SL(2,\Complex)$ quantum Hilbert space with $\Hilb$ depends on a choice of complex structure $\tau$ on $\Sigma$.
Such a choice is not invariant under the action of $\Gamma$, indeed $\gamma\in\Gamma$ acts on $\Half$ via Mobius transformations.
Nevertheless we can compute the action of $\tilde S$ and $\tilde T$ on $\torus\oplus\torus$ as follows
\begin{align}\label{eq:STaction}
\tilde S \psi (\theta_1,\theta_2) &= \psi(\theta_2,-\theta_1)\\
\label{eq:STaction2}
\tilde T \psi (\theta_1,\theta_2) &= \psi (\theta_1,\theta_1+\theta_2).
\end{align}
We remark that these operators commute with the action of the Weyl group and are equivariant with respect to the action of $\coroot\times\coroot$ (see proof of Lemma A.2 in \cite{Je}).

  \section{Weil-Gel'fand-Zak transform on lattices.} \label{sc:WGZ}
Let $\left(E,\rprod{\cdot}{\cdot}\right)$ by a fixed euclidean space of dimension $\dim E = n<\infty$. Let $\Lambda\subset E$ be a full-rank lattice. Denote $\Torus_E$ the quotient $E/\Lambda$.
Let $L\rightarrow \Torus_E\oplus \Torus_E$ be the line bundle defined by the following multipliers
\begin{align}
e_{\lambda_1\oplus\lambda_2}(x\oplus y) =(-1)^{\rprod{\lambda_1}{\lambda_2}} e^{-\pi i( \rprod{x}{\lambda_2} - \rprod{\lambda_1, y})}
\end{align}
Define the following notation for  Fourier kernels
\begin{align}
\fker{x}{y} \equiv e^{ 2\pi i\rprod{x}{y}}
\end{align} 
and define the dual $\Lambda^*\subset E$ as the lattice
\begin{align}
\Lambda^*\equiv \left\{\gamma\in E\text{ : } \fker{\lambda}{\gamma} = 1\text{, }\forall\, \gamma\in\Lambda \right\}.
\end{align}
We say that $\Lambda$ is \emph{integral} if $\Lambda\subseteq \Lambda^*$. In such case $\Zcal \equiv \Lambda^*/\Lambda$ is a well defined finite abelian group.\\

Suppose now $\Lambda$ is integral.
The space $\SqInt(E)\otimes \setC^{\Zcal}$ has an inner product defined as follows. Let $f_\gamma,g_\gamma\in\SqInt(E)$  for $\gamma\in\Zcal$ be square integrable functions on $E$, for 
\begin{align*}
\mathsf{f} = \left(f_\gamma\right)_{\gamma\in\Zcal}\text{, and }\mathsf{g} = \left(g_\gamma\right)_{\gamma\in\Zcal}\,\in\SqInt(E)\otimes \setC^{\Zcal}
\end{align*}

 we define
\begin{equation}\label{eq:L2product}
\inn{\mathsf f}{\mathsf g} = \sum_{\gamma\in\Zcal} \int_{E}f_\gamma\conj{g_\gamma}\dvol,
\end{equation}
where the volume form $\dvol$ is specified by the inner product $\rprod{\cdot}{\cdot}$. We write 
$$\Vol{\Lambda} = \int_{F_{\Lambda}}\dvol$$
for the volume of a fundamental domain $F_{\Lambda}\subset E$ for the action of $\Lambda$.\\
On the space of sections $\SqInt(\Torus_E\times\Torus_E, L)$ we consider the inner product 
\begin{align}\label{eq:L2prodTori}
\inn{\psi}{\phi} = \frac{1}{\vol{\Lambda}}\int_{F_{\Lambda}\times F_{\Lambda}}\psi(\theta_1,\theta_2)\conj{\phi(\theta_1,\theta_2)}\dvol(\theta_1,\theta_2)
\end{align}
where the measure is the product measure on $E\times E$.\\

Let $S(E)$ be the space of Scwartz functions on $E$.
The following is the {\em Weil-Gel'fand-Zak transform on the lattice $\Lambda$}.
\begin{theorem}\label{pr:WGZ}
	\label{prop:1}
	We have an isomorphism
	\begin{align}
	\nonumber
	Z\colon S(E)\otimes \setC^\Zcal \longrightarrow
	C^\infty(\Torus_E\times\Torus_E, L)
	\end{align}
	given by
	\begin{align}
	Z(\mathsf f)(\theta_1,\theta_2) &=\frac{1}{\sqrt{\abs{\Zcal}}} \fker{-{\theta_1}/{2}}{\theta_2}
	\sum_{\gamma\in\Zcal}\sum_{\lambda\in\Lambda^*} f_\gamma(\theta_1 + \lambda) \fker{-\lambda}{\theta_2}\fker{-\lambda}{\gamma},
	\end{align}
	 with inverse
	\begin{align}
	\nonumber
	 Z^{-1} (s) (\theta, \gamma) &= \frac{1}{\sqrt{\abs{\Zcal}}} \sum_{\hat{\gamma}\in\Zcal} \fker{\gamma}{\hat{\gamma}} \frac{1}{\vol{\Lambda}}\int_{F_\Lambda} s(\theta - \hat{\gamma}, \tilde\theta) \fker{\theta+\hat{\gamma}}{\tilde\theta/2} \dvol(\tilde\theta)\text{.}
	\end{align}
	It satisfies the following unitarity property
	\begin{align}
	\nonumber
	\inn{Z(\mathsf f)}{Z( \mathsf{g})} = \inn{\mathsf f}{ \mathsf g}
	\end{align}
	thus $Z$ extends to an isometry between $\SqInt(E)\otimes \setC^{\Zcal}$ and $\SqInt(\Torus_E\times\Torus_E, L)$.
\end{theorem}
\proof
We write $\diff\theta$ instead of $\dvol(\theta)$ in the computations.
\begin{align*}
Z\left(Z^{-1} (s)\right) (\theta_1,\theta_2)
&= \frac{1}{\sqrt{\abs{\Zcal}}} \fker{-{\theta_1}/{2}}{\theta_2} \sum_{\gamma\in\Zcal}\sum_{\lambda\in\Lambda^*} Z^{-1} (s) (\theta_1+\lambda, \gamma)\fker{-\lambda}{\theta_2}\fker{-\lambda}{\gamma}\\
&= \frac{1}{{\abs{\Zcal}}} \fker{-{\theta_1}/{2}}{\theta_2} \sum_{\hat{\gamma}\in\Zcal}\sum_{\lambda\in\Lambda^*}  \frac{1}{\vol{\Lambda}}
\fker{-\lambda}{\theta_2}
\sum_{{\gamma}\in\Zcal}\fker{\gamma}{\hat \gamma- \lambda}\\
&\qquad\times  \int_{F_\Lambda} 
s(\theta_1+\lambda - \hat\gamma, \tilde{\theta})\fker{\theta_1 +\lambda+\hat{\gamma}}{\tilde{\theta}/2}\diff \tilde{\theta}\\
&=\fker{-{\theta_1}/{2}}{\theta_2} \sum_{\tilde{\gamma}\in\Zcal}\sum_{\tilde\lambda\in\Lambda}\frac{1}{\vol{\Lambda}}
\int_{F_\Lambda} 
s(\theta_1+\tilde\lambda, \tilde{\theta})\fker{\theta_1 +\tilde\lambda+2\tilde{\gamma}}{\tilde{\theta}/2}\diff \tilde{\theta}\\
&\qquad\times \fker{-\tilde{\lambda} -\tilde{\gamma}}{\theta_2}
\\
&=\fker{-{\theta_1}/{2}}{\theta_2} \sum_{\tilde{\gamma}\in\Zcal}\sum_{\tilde\lambda\in\Lambda}\frac{1}{\vol{\Lambda}}
\int_{F_\Lambda} 
s(\theta_1, \tilde{\theta})\fker{\frac{\theta_1}{2} +\tilde\lambda+\tilde{\gamma}}{\tilde{\theta}}\diff \tilde{\theta}\\
&\qquad\times \fker{-\tilde{\lambda} -\tilde{\gamma}}{\theta_2}\\
&= s(\theta_1, \theta_2),
\end{align*}
where the last equation is given by the Fourier series properties of the $\Lambda$-periodic function $s(\theta_1,\tilde{\theta})\fker{\theta_1/2}{\tilde \theta}$ in the variable $\tilde{\theta}$. 

To establish the unitarity of the transform, we calculate
\begin{align*}
\inn{Z(\mathsf f)}{Z( \mathsf{g})} &=  \frac{1}{\vol{\Lambda}}\int_{F_{\Lambda}\times F_{\Lambda}}
Z(\opf)(\theta_1,\theta_2)\conj{Z(\opg)(\theta_1,\theta_2)} \diff\theta_1\diff\theta_2\\
&=\frac{1}{\vol{\Lambda}}\frac{1}{{\abs{\Zcal}}}\sum_{\gamma,\hat{\gamma}\in\Zcal}\sum_{\lambda,\hat\lambda\in\Lambda^*}\fker{-\lambda}{\gamma}\fker{\hat\lambda}{\hat\gamma}
\int_{F_{\Lambda}}
\opf(\theta_1+\lambda,\gamma)\conj{\opg(\theta_1+\hat{\lambda},\hat{\gamma})} \diff\theta_1\\
&\qquad\times
\int_{F_{\Lambda}}\fker{\hat{\lambda}-\lambda}{\theta_2}\diff\theta_2\\
&=\frac{1}{{\abs{\Zcal}}}\sum_{\gamma,\hat{\gamma}\in\Zcal}\sum_{\lambda\in\Lambda^*}\fker{-\lambda}{\gamma}\fker{\lambda}{\hat\gamma}
\int_{F_{\Lambda}}
\opf(\theta_1+\lambda,\gamma)\conj{\opg(\theta_1+{\lambda},\hat{\gamma})} \diff\theta_1\\
&=\frac{1}{{\abs{\Zcal}}}\sum_{\gamma,\tilde{\gamma},\hat{\gamma}\in\Zcal}\sum_{\tilde\lambda,\in\Lambda}\fker{-\tilde{\gamma}}{\gamma}\fker{\tilde{\gamma}}{\hat\gamma}
\int_{F_{\Lambda}}
\opf(\theta_1+\tilde\lambda + \tilde{\gamma},\gamma)\conj{\opg(\theta_1+\tilde{\lambda} +\tilde{\gamma},\hat{\gamma})} \diff\theta_1\\
&=\sum_{\tilde{\gamma}\in\Zcal}
\int_{E}
\frac{1}{\sqrt{\abs{\Zcal}}}\sum_{{\gamma}\in\Zcal}\opf(\theta_1,\gamma)\fker{-\tilde{\gamma}}{\gamma}
\frac{1}{\sqrt{\abs{\Zcal}}}\sum_{{\hat\gamma}\in\Zcal}\conj{\opg(\theta_1,\hat{\gamma})} \fker{\tilde{\gamma}}{\hat\gamma}\diff\theta_1\\
&=\sum_{\tilde{\gamma}\in\Zcal}
\int_{E}\opf(\theta_1,\tilde \gamma)\conj{\opg(\theta_1,\tilde{\gamma})}\diff\theta_1
\end{align*}
\endproof

We define the following auxiliary operators on $\SqInt(E) \otimes \Complex^{\Zcal}$
\begin{align}
\label{eq:FiniteOps}
&\opF_{\Zcal}(\mathsf f)(\theta,\gamma) = \frac{1}{\sqrt{\abs{\Zcal}}} \sum_{\hat{\gamma}\in\Zcal} f_\gamma(\theta) \fker{\gamma}{\hat{\gamma}}
&\opG_{\Zcal}(\mathsf f)(\theta,\gamma) = \fker{\gamma/2}{\gamma} \mathsf f(\theta, \gamma)\\
\label{eq:ContOps}
&\mathcal{F}_{E}(\mathsf f)(\theta,\gamma) = \int_{E} f_\gamma(\tilde{\theta})\fker{\theta}{\tilde{\theta}}\dvol(\theta)
&\mathcal{G}_{E} (\mathsf f)(\theta,\gamma) = \fker{\theta/2}{\theta}(\mathsf f)(\theta,\gamma)\text{.}
\end{align}
Let $\tilde S$ and $\tilde T$ be two operators acting on $\SqInt(\Torus_E\times\Torus_E, L)$ as expressed in (\ref{eq:STaction}-\ref{eq:STaction2}). Define the following operators acting on $\SqInt(E)\otimes \setC^{\Zcal}$
\begin{align}
&\hat S = \left(Z\circ\opF_{\Zcal}\right)^{-1} \circ\tilde S \circ \left(Z\circ\opF_{\Zcal}\right),
&\hat{T} =  \left(Z\circ\opF_{\Zcal}\right)^{-1} \circ\tilde T \circ \left(Z\circ\opF_{\Zcal}\right).
\end{align} 
They can be explicitly computed as
\begin{proposition}\label{pr:pqact}
\begin{align}
\hat S (\mathsf f) (\theta, \gamma) &=\opF^{-1}_{\Zcal} \circ\mathcal{F}_{E} (\mathsf f) (\theta,\gamma),\\
\hat T (\mathsf f) (\theta, \gamma) &= \opG_{\Zcal} \circ \mathcal{G}^{-1}_{E} (\mathsf f) (\theta, \gamma).
\end{align}
\end{proposition}
\proof 
\begin{align*}
Z^{-1}\circ \tilde S\circ Z(\opf)(\theta,\gamma) &= 
\frac{1}{\sqrt{\abs{\Zcal}}} \sum_{\hat{\gamma}\in\Zcal} \fker{\gamma}{\hat{\gamma}} \frac{1}{\vol{\Lambda}}\int_{F_\Lambda} S(Z(\opf))(\theta - \hat{\gamma}, \tilde\theta) \fker{\theta+\hat{\gamma}}{\tilde\theta/2} \diff\tilde\theta\\
&= \frac{1}{\sqrt{\abs{\Zcal}}} \sum_{\hat{\gamma}\in\Zcal} \fker{\gamma}{\hat{\gamma}} \frac{1}{\vol{\Lambda}}\int_{F_\Lambda} Z(\opf)(\tilde\theta , -\theta + \hat{\gamma}) \fker{\theta+\hat{\gamma}}{\tilde\theta/2} \diff\tilde\theta\\
&= \frac{1}{\vol{\Lambda}}\frac{1}{{\abs{\Zcal}}} \sum_{\tilde\gamma\in\Zcal}\sum_{\lambda\in\Lambda^*}
\int_{F_\Lambda}\opf(\tilde{\theta}+\lambda, \tilde{\gamma})\fker{\theta}{\tilde{\theta}}\fker{\lambda}{\theta - \tilde{\gamma}}\diff \tilde{\theta}\\
&\qquad\times\sum_{\hat\gamma\in\Zcal}\fker{\hat{\gamma}}{\gamma -\lambda}\\
&=\frac{1}{\vol{\Lambda}} \sum_{\tilde\gamma\in\Zcal}\sum_{\tilde\lambda\in\Lambda}
\int_{F_\Lambda}\opf(\tilde{\theta}+\tilde\lambda + \gamma, \tilde{\gamma})\fker{\theta}{\tilde{\theta} + \tilde\lambda + \gamma}\diff \tilde{\theta} \fker{\gamma}{- \tilde{\gamma}}\\
&=\frac{1}{\vol{\Lambda}} \sum_{\tilde\gamma\in\Zcal}
\int_{E}\opf(\tilde{\theta}, \tilde{\gamma})\fker{\theta}{\tilde{\theta} }\diff \tilde{\theta} \fker{\gamma}{- \tilde{\gamma}}.
\end{align*}
Noticing that $\abs{\Zcal} = \left(\vol{\Lambda}\right)^2$ we get the first equation.
\begin{align*}
Z^{-1}\circ \tilde T\circ Z(\opf)(\theta,\gamma) &= 
\frac{1}{\sqrt{\abs{\Zcal}}} \sum_{\tilde{\gamma}\in\Zcal} \fker{\gamma}{\tilde{\gamma}} \frac{1}{\vol{\Lambda}}\int_{F_\Lambda} T(Z(\opf))(\theta - \tilde{\gamma}, \tilde\theta) \fker{\theta+\tilde{\gamma}}{\tilde\theta/2} \diff\tilde\theta\\
&= 
\frac{1}{\sqrt{\abs{\Zcal}}} \sum_{\tilde{\gamma}\in\Zcal} \fker{\gamma}{\tilde{\gamma}} \frac{1}{\vol{\Lambda}}\int_{F_\Lambda} Z(\opf)(\theta - \tilde{\gamma}, \tilde\theta + \theta -\tilde{\gamma}) \fker{\theta+\tilde{\gamma}}{\tilde\theta/2} \diff\tilde\theta\\
&= 
\frac{1}{{\abs{\Zcal}}} \frac{1}{\vol{\Lambda}} \sum_{\hat\gamma, \tilde{\gamma}\in\Zcal} \fker{\gamma}{\tilde{\gamma}} 
\int_{F_\Lambda} \sum_{\lambda\in\Lambda^*} \opf(\theta - \tilde{\gamma} + \lambda, \hat{\gamma}) 
\fker{-\frac{\theta - \tilde{\gamma}}{2}}{\theta -\tilde{\gamma}+\tilde{\theta}}\\
&\qquad\times \fker{-\lambda}{\theta-\tilde{\gamma}+\tilde{\theta}}\fker{-\lambda}{\hat{\gamma}}
 \fker{\theta+\tilde{\gamma}}{\tilde\theta/2} \diff\tilde\theta\\
&= 
\frac{1}{{\abs{\Zcal}}}  \sum_{\hat\gamma, \tilde{\gamma}\in\Zcal} \fker{\gamma}{\tilde{\gamma}} 
 \sum_{\lambda\in\Lambda^*} \opf(\theta - \tilde{\gamma} + \lambda, \hat{\gamma}) 
\fker{-\frac{\theta - \tilde{\gamma}}{2}}{\theta -\tilde{\gamma}}\\
&\qquad\times \fker{-\lambda}{\theta-\tilde{\gamma}}\fker{-\lambda}{\hat{\gamma}}
\frac{1}{\vol{\Lambda}}\int_{F_\Lambda}\fker{\tilde{\theta}}{\tilde{\gamma }- \lambda}\diff\tilde\theta\\
&= 
\frac{1}{{\abs{\Zcal}}}  \sum_{\hat\gamma, \tilde{\gamma}\in\Zcal} \fker{\gamma}{\tilde{\gamma}} 
\opf(\theta, \hat{\gamma}) 
\fker{-\frac{\theta - \tilde{\gamma}}{2}}{\theta -\tilde{\gamma}} \fker{-\tilde{\gamma}}{\theta-\tilde{\gamma}}\fker{-\tilde{\gamma}}{\hat{\gamma}}\\
&= 
\fker{-\theta/2}{\theta}\frac{1}{\sqrt{\abs{\Zcal}}}  \sum_{ \tilde{\gamma}\in\Zcal} \fker{\gamma}{\tilde{\gamma}} \fker{\frac{\tilde{\gamma}}{2}}{\tilde{\gamma}}
\frac{1}{\sqrt{\abs{\Zcal}}}\sum_{ \hat{\gamma}\in\Zcal}
\opf(\theta, \hat{\gamma}) \fker{-\tilde{\gamma}}{\hat{\gamma}}\\
&= \mathcal{G}_E^{-1}\circ \opF_{\Zcal}\circ\opG_{\Zcal}\circ\opF_{\Zcal}^{-1}(\opf)(\theta, \gamma)
\end{align*}

\endproof
It is evident that both $\hat{S}$ and $\hat{T}$ have a (unique up to a scalar) tensor product decomposition corresponding to $\SqInt(E)\otimes\Complex^{\Zcal}$. We write this precisely in the next proposition

\begin{proposition}\label{pr:pqactnomr}
Let $j$ and $\omega$ be complex numbers satisfying
\begin{align}
&j^4 =1,
&\omega^3 = i^{n/2}j^{-1}.
\end{align}
The operators $\hat S$ and $\hat T$ generate a representation of $SL(2,\mathbb Z)$ as endomorphisms of $\SqInt(E)\otimes\Complex^{\Zcal}$. For any possible choice of $j$ and $\omega$ there is a tensor product decomposition
\begin{align}
&\hat S = \hat{S}^\p\otimes\hat{S}^\pp
&\hat{T} = \hat{T}^\p\otimes\hat{T}^\pp\\
&\hat{S}^\p = j\mathcal{F}_E
&\hat{T}^\p = \omega\mathcal G^{-1}_{E}\\
&\hat{S}^\pp = j^{-1} \opF_{\Zcal}
&\hat{T}^\pp = \omega^{-1} \opG_{\Zcal}
\end{align}
where $\hat{S}^\p$ and $\hat{T}^\p$ (resp. $\hat{S}^\pp$ and $\hat{T}^\pp$ ) generates a representation on $\SqInt(E)$ (resp. $\Complex^{\Zcal}$)
\end{proposition}
\proof
The proof is just a direct verification of the relations of $SL(2,\Integer)$.
\endproof

We now apply these computations to the genus $1$ Chern--Simons theory presented in the previous section (see also Appendix \ref{ap:Lie} for further notation). The Euclidean space is $\left(\torus, \krprod{\cdot}{\cdot}\right)$ while the line bundle is $\pi^*\Line^k\rightarrow \Torus\times\Torus$. The integral lattice we consider is $\coroot$ and its dual will be $\kweights$. Their quotient is the finite abelian group $\kZcal$.  The operator $Z\circ\opF_{\kZcal}$ defines an isometry between $\widetilde{\Hilb}^{(k)}$ and $\SqInt(\torus)\otimes\Complex^{\kZcal}$.\\

In this way, Proposition \ref{pr:pqactnomr} provides an explicit description of a lift of the action of the mapping class group $\Gamma$ to $\widetilde{\Hilb}^{(k)}$.\\
Consider the following action of the Weyl group $W$ on $\SqInt(\torus)\otimes\Complex^{\kZcal}$
\begin{align}\label{eq:WhatAct}
\hat{w}\cdot \opf(\theta, \gamma) \equiv f(w(\theta), w(\gamma))\text{, }\forall\, w\in W.
\end{align}
The operator $(Z\circ \opF_{\kZcal})$ is equivariant with respect to this action and the one on $\widetilde{\Hilb}^{(k)}$. Therefore it defines an isometry
\begin{align}
Z\circ \opF_{\kZcal} \colon \Hilb^{(k)} \longrightarrow \left(S(\torus)\otimes\Complex^{\kZcal}\right)^W.
\end{align}
The action \eqref{eq:WhatAct} behaves well with respect to the tensor product decomposition $\SqInt(\torus)\otimes\Complex^{\kZcal}$. Therefore it can be written as the tensor product $\hat{w} = \hat{w}^\p\otimes\hat{w}^\pp$ of the two actions 
\begin{align}
&\hat{w}^\p(f) (\theta) =  f(w(\theta))
&\hat{w}^\pp (x)(\gamma) = x(w(\gamma)).
\end{align}
On a $W$--invariant vector $f\otimes x$,  i.e. $\hat{w}\cdot f\otimes x = f\otimes x$, there will exist $\lambda_w\in\Complex^*$ such that
\begin{align}
&\hat{w}^\p(f) (\theta) = \lambda_w f(\theta)
&\hat{w}^\pp (x)(\gamma) = \lambda_w^{-1} x(\gamma).
\end{align}
Since $W$ is generated by elements of order $2$, $\lambda_{w} \in \{1,-1\}$. In fact $\lambda_w$ has to be either the identical character $1$ or the alternating character $\det(w)$. 
In the end, a vector $f\otimes x$ is $W$ invariant if and only if  $f$ and $x$ are both invariant or both anti-invariant for the actions of $\hat{w}^\p$ and $\hat{w}^\pp$ respectively. With the notation from Appendix \ref{ap:Lie} we can write
\begin{align}
\nonumber
\Hilb^{(k)}&\simeq \left(\SqInt(\torus)\otimes\Complex^{\kZcal}\right)^W \\
\label{eq:WeylDec}
&= \left(\SqInt_{inv}(\torus)\otimes \overline{C}_{(k)}\right)\bigoplus\left(\SqInt_{anti}(\torus)\otimes {C}_{(k)}\right)
\end{align}

Moreover, all the operators involved in the description of the action of the mapping class group from proposition \ref{pr:pqactnomr} are invariant with respect to the action \eqref{eq:WhatAct}, so they preserves the whole decomposition \eqref{eq:WeylDec}.\\

On tensors $f\otimes x$, and $g\otimes y\in \SqInt(\torus)\otimes\Complex^{\kZcal}$ the inner product $\inn{\cdot}{\cdot}$ from \eqref{eq:L2product} factor into two inner products
\begin{align}\label{eq:inners}
&\inn{f}{g}_{\torus} \equiv \int_{\torus}f(\theta)\conj{g(\theta)}\dvol_k(\theta),
&\inn{x}{y}_{\kZcal} \equiv \sum_{\gamma\in\kZcal}x(\gamma)\conj{y(\gamma)}
\end{align}
We can use the second one, together with the basis from \eqref{eq:InvFinBasis} and \eqref{eq:AntiFinBasis}, to  compute the matrix elements of the action of SL$(2,\Integer)$ expressed in Proposition \ref{pr:pqactnomr}, on the spaces $\overline{C}_{(k)}$ and $C_{(k)}$
\begin{lemma}\label{lm:finiteOps}
\begin{align*}
&\inn{e_\alpha}{j^{-1}\opF_{\kZcal} e_\beta}_{\kZcal} = \frac{j^{-1}}{\sqrt{\abs{\kZcal}} }
\sum_{w\in W} \fker{w\alpha}{\beta},
&\inn{e_\alpha}{\omega^{-1}\opG_{\kZcal}^{-1} e_\beta}_{\kZcal} = \omega^{-1} \fker{-\frac{\alpha}{2}}{\alpha}\delta_{\alpha,\beta},\\
 & \inn{\tilde e_\alpha}{\opF_{\kZcal} \tilde e_\beta}_{\kZcal} = \frac{j^{-1}}{\sqrt{\abs{\kZcal}} }
 \sum_{w\in W} \det(w)\fker{w\alpha}{\beta},
 & \inn{\tilde e_\alpha}{\omega^{-1}\opG_{\kZcal}^{-1} \tilde e_\beta}_{\kZcal} = \omega^{-1}\fker{-\frac{\alpha}{2}}{\alpha}\delta_{\alpha,\beta},\\
&\inn{\tilde e_\alpha}{\omega^{-1}\opG_{\kZcal}^{-1} e_\beta}_{\kZcal} = \inn{\tilde e_\alpha}{\opF_{\kZcal}  e_\beta}_{\kZcal} = 0\text{.}
\end{align*}
\end{lemma}
\proof
This is a direct verification from the definitions.
\endproof
\emph{Proof of Theorem \ref{th:pqr}}
Lemma \ref{lm:finiteOps} together with Proposition \ref{pr:pqactnomr} provides a proof for the decomposition and explicit formulas of $\hat{\rho}_k$ in Theorem \ref{th:qr}. 
$\hfill \square$\\

Before proceeding with a review of the genus $1$ Hitchin-Witten connection, we describe the pre-quantum  connection $\nabla$ from equation (\ref{eq:pqConnection}), as an explicit differential operators acting on $\SqInt(\torus)\otimes\Complex^{\kZcal}$.
Let namely $\hat{\nabla}_X = (Z\circ \opF_{\kZcal})^{-1} \circ\nabla_X \circ (Z\circ \opF_{\kZcal})$ for any $X\in\Smooth(\Torus\times\Torus, T(\Torus\times\Torus))$.
Define the following operator
\begin{align}
\nonumber
D_{\sigma, j} \colon S(\torus)\otimes\Complex^{\kZcal} &\longrightarrow S(\torus)\otimes\Complex^{\kZcal}\\
\mathsf f(\theta,\gamma)&\mapsto \diff f_\gamma\left[\opb_j\right](\theta) +2\pi i\,\conj\sigma^{-1}\krprod{{\opb_j}}{\theta}f_\gamma(\theta)\label{eq:DsigmaDef}
\end{align}
Notice that theses operators depend on the basis $\{\opb_j\}_{j}$ of $\torus$.\\

From now on, given a connection $\nabla$ and a local coordinate function $u$ we will write $\nabla_u$ in place of $\nabla_{\frac{\partial}{\partial u}}$. Similarly we will write $\partial_u$ in place of $\frac{\partial}{\partial u}$.
We have that
\begin{lemma}\label{lm:hatNabla}
	 Let $u_j$ and $v_j$, $\sigma\in\Half$, $z_j$ and $\conj{z}_j$, for $j=1,\dots,n$ be defined as in \eqref{eq:thetacoor} and \eqref{eq:holocoord}.
	 We have that
	\begin{align}
	&\hat{\nabla}_{{{u_j}}} f_\gamma (\theta) = \diff f_\gamma\left[\opb_j\right](\theta)
	&\hat{\nabla}_{{{v_j}}} f_j(\theta) = -2\pi i\krprod{{\opb_j}}{\theta} f_\gamma(\theta)\\
	&\hat{\nabla}_{{{z_j}}} f_j(\theta) =\frac{\conj \sigma}{\conj \sigma - \sigma} D_{\sigma, j} f_\gamma(\theta)
	&\hat{\nabla}_{{\conj{z}_j}} f_\gamma(\theta) = \frac{\sigma}{\sigma - \conj{\sigma}}D_{\conj{\sigma}, j} f_\gamma(\theta)
	\end{align}
\end{lemma}
\proof 
This is a simple verification using the explicit formulas for $\nabla$ and $Z$.
\endproof
We remark that the pre-quantum connection $\hat{\nabla}$ decomposes under the tensor decomposition $\SqInt(\torus)\otimes\Complex^{\kZcal}$ and it is trivial in the second factor. In particular the Hitchin--Witten parallel transport does not affect the seond factor.

\section{The Hitchin-Witten Connection in genus $1$}\label{sc:HWC}
The quantum Hilbert space for Chern-Simons theory with gauge group $G_\Complex$ can be identified with the pre-quantum Hilbert space for Chern Simons theory with gauge group $G$, as explained in \cite{W3}. 
This identification, however, depends on a choice of a complex structure on the surface $\Sigma$.
In genus $1$ this means that there is a trivial bundle
\begin{equation}\label{eq:HWbundle}
\Half\times\Hilb^{(k)} \longrightarrow \Half
\end{equation}
and the identification between the fibers is obtained through parallel transport with the Hitchin-Witten connection $\boldnabla$ defined in \cite{W3} for genus $g$ surfaces and in \cite{AG2} in a more general setting.
This connection is in this case flat, so the identifications between different quantum spaces will be well defined.
The explicit formula for genus $1$ is
\begin{align}
\nonumber
\boldnabla_{{\partial_\sigma}} &= \partial_\sigma +
\frac{1}{2t}\Delta_{G} \\
\nonumber
\boldnabla_{{\partial_{\conj\sigma}}} &=\partial_{\conj\sigma} -
\frac{1}{2\conj t}\Delta_{\overline{G}}
\end{align}
where $\sigma$ is the holomorphic coordinate for $\Half$ and $t\in\Complex^*$ with $t = k +is$.
The operators  $\Delta_{G}$ and $\Delta_{\overline{G}}$ are second order differential operators on the line bundle $\Line^k$ defined in \cite{AnB} and \cite{AG2}. Explicitly, for genus $1$, we have that
\begin{align*}
&\Delta_{G} = \frac{i}{\pi}\sum_{p,q=1}^nC^{p,q}_{(1)} \nabla_{z_p}\nabla_{z_q}
&\Delta_{\conj G} = -\frac{i}{\pi} \sum_{p,q=1}^n C^{p,q}_{(1)} \nabla_{\conj z_p}\nabla_{\conj z_q}
\end{align*}
where $C^{p,q}_{(k)}$ is the $p,q$ entry of the inverse matrix of $C_{l,m}^{(k)}= \krprod{{\opb_l}}{{\opb_m}}$.
\begin{proposition}
	Suppose $s\in\Real$, $\phi$, $\psi\in \Hilb^{(k)}$ and $\inn{\cdot}{ \cdot}: \Hilb^{(k)}\times\Hilb^{(k)}\rightarrow \Complex$ as in (\ref{eq:hermitianDef}). We then have that
	\begin{equation*} 
	\diff\inn{\phi}{\psi} = \inn{\boldnabla\phi}{ \psi} + \inn{\phi}{ \boldnabla\psi}.
	\end{equation*}
\end{proposition}
\proof
Partial integration together with the hypothesis $s\in\Real$ gives
\begin{equation}
\inn{\phi}{\frac{1}{2t}\Delta_G \psi} =
  \inn{\frac{1}{2\conj t}\Delta_{\conj G} \phi}{ \psi} .
\end{equation}
A direct computation of the sum $\inn{\boldnabla\phi}{ \psi} + \inn{\phi}{ \boldnabla\psi}$ using the equation above gives the result.
\endproof 

Let us define the Laplace operator operator 
\begin{equation}\label{eq:Laplacian}
\Delta_{\sigma} = \frac{i}{2\pi}\left(\sigma - \conj\sigma \right)\sum_{p,q= 1}^{n}\,C^{p,q}_{(1)} \left(\nabla_{z_p}\nabla_{\conj z_q}
+ \nabla_{\conj z_p}\nabla_{z_q}\right)
\end{equation}
which clearly dependent on the complex structure determined by $\sigma\in {\mathbb H}$.
\begin{lemma}\label{lm:commRules}
On $\Hilb^{(k)}\times\Half$, the following hold true
\begin{align}
&\comm{\nabla_{z_p}\text{,} \, \nabla_{\conj z_q}} =\frac{2\pi i C^{(k)}_{p,q}}{\sigma -\conj{\sigma}}\\
&\comm{\partial_\sigma\text{,} \nabla_{z_p}} = -\comm{\partial_\sigma\text{,} \nabla_{\conj z_p}} 
= \frac{-1}{\sigma - \conj{\sigma}} \nabla_{z_p}\text{, }\\
&\comm{\partial_{\conj\sigma}\text{,} \nabla_{\conj z_p}} = -\comm{\partial_{\conj\sigma}\text{,} \nabla_{z_p}} 
= \frac{-1}{\conj\sigma - {\sigma}} \nabla_{\conj z_p}\\
&\comm{\partial_{\sigma}\text{,} \, \Delta_\sigma} = \Delta_{G}\text{,  }\qquad\qquad\qquad\qquad \quad\,\,\,\,
\comm{\partial_{\conj\sigma}\text{,} \, \Delta_\sigma} = \Delta_{\conj G}\\
&\comm{\Delta_\sigma\text{,} \,  \Delta_{G}} = 4k  \Delta_{G}\text{,    }\,\,\qquad\qquad \qquad\qquad
\comm{\Delta_\sigma\text{,} \,  \Delta_{\conj G}} = -4k  \Delta_{\conj G}
\end{align}
\end{lemma}
\proof
The first commutator follows from the explicit curvature of $\nabla$. The others follows from iterations of the first one or from the explicit dependence of $z_p$ on $\sigma$.
\endproof

\begin{lemma}\label{lm:algebraComm}
	Let $x$, $y$ and $z$ be three elements of a Lie algebra satisfying the relations $\comm{x\text{,}y} = z$ and $\comm{y,z} = az$ for a central $a$. Suppose that the formal exponential $e^y = \sum_{k\geq 0} \frac{y^k}{k!}$ is well defined. Then we have
	\begin{equation*}
	\comm{x,e^y} = \frac{z}{a} \left( e^{y+a} - e^y\right)
	\end{equation*} 
\end{lemma}
As was noted by Witten \cite{W2}, Lemma \ref{lm:algebraComm} together with Lemma \ref{lm:commRules} gives us the following conjugation rule
\begin{equation}
e^{-r\Delta_{\sigma}} \circ \diff \circ e^{r\Delta_{\sigma}} = \boldnabla
\end{equation}
where $r$ is chosen so that 
\begin{equation}\label{eq:rDef}
e^{-4kr} = -\frac{k- is}{k+is}\text{.} 
\end{equation}
In particular, equation (\ref{eq:HWconj}) implies 
\begin{proposition}\label{pr:parallelSections}
For every $\psi\in\Hilb^{(k)}$ independent on the complex structure of $M$, the section
$$ e^{-r\Delta}\psi$$
of the vector bundle (\ref{eq:HWbundle}) is parallel with respect to $\boldnabla$.
\end{proposition}

Given $\gamma\in\Gamma$, its pre-quantum action on $\Hilb^{(k)}$ was defined in (\ref{eq:STaction}), however when we look at the action on the whole bundle $\Hilb^{(k)}\times \Half\rightarrow\Half$, $\gamma$ acts on $\Half$ by  $\gamma_*$ as pull-back via the Mobius transformation $\gamma^{-1}$. We will then need to compose $\tilde \gamma$ with the parallel transport $\mathcal P_{ \sigma, \gamma_*\sigma}$ of the pre-quantum action with the Hitchin-Witten connection from $\gamma_*\sigma$ back to $\sigma$.
By the results in  Proposition \ref{pr:parallelSections}, we have that
\begin{equation}
\mathcal P_{\sigma_0, \sigma_1} \psi_{\sigma_1} = e^{-r\Delta_{\sigma_0}} e^{r\Delta_{\sigma_1}}\psi_{\sigma_1},
\end{equation}
so
\begin{equation}
\eta_{t}(\varphi) = e^{-r\Delta_{\sigma}} e^{r\Delta_{\varphi_*\sigma}}\circ \rho_k(\varphi).
\end{equation}
If $w = x+\tau y$ are holomorphic coordinates for the surface $\Sigma$, then we have that
\begin{align*}
&S_* \tau = -\frac{1}{\tau}  &T_*\tau = \tau -1\text{,}
\end{align*}
and recalling that the holomorphic coordinate $z_p= u_p + \sigma v_p$ on $M$ are related with $\tau$ by $\sigma = -\frac 1\tau$, we get that
\begin{align*}
&S_* \sigma = -\frac{1}{\sigma}  &T_*\sigma = \frac{\sigma}{1+\sigma}\text{.}
\end{align*}

However the Laplace operator is itself is mapping class group invariant, as one can directly verify that 
\begin{align}
&\tilde{T}\circ\Delta_{\sigma} = \Delta_{T_*\sigma}\circ \tilde T,
&\tilde{S}\circ\Delta_{\sigma}= \Delta_{S_*\sigma}\circ \tilde{S}.
\end{align}
For any $\varphi\in\Gamma$ we obtain the equation
\begin{equation}\label{eq:qr}
e^{-r\Delta_{\sigma}} e^{r\Delta_{\varphi_*\sigma}}\circ \tilde{\varphi} =
 e^{-r\Delta_{\sigma}}\circ \tilde{\varphi} \circ e^{r\Delta_{\sigma}}
\end{equation}
and this proves Proposition \ref{pr:qr}. \\

Let $\hat{\Delta}_\sigma \equiv \left(Z\circ\opF_{\Zcal}\right)^{-1}\circ\Delta_\sigma\circ \left(Z\circ\opF_{\Zcal}\right)$. We remark that, as the connection $\hat{\nabla}$ acts as the identity on the second factor of the tensor product $\SqInt(\torus)\otimes\Complex^{\kZcal}$, we can consider $\hat{\Delta}_\sigma$ as an operator on $\SqInt(\torus)$ alone.
\begin{proposition}\label{pr:DeltaSpectra}
	We have that
	\begin{equation}\label{eq:hatDelta2}
	\hat{\Delta}_{\sigma} =nk+ \frac{i}{\pi}\sum_{p,q=1}^{n} \frac{\sigma\conj{\sigma}}{\conj{\sigma}-\sigma} C^{p,q}_{(1)}D_{\sigma,p}\,D_{\conj{\sigma}, q}.
	\end{equation}
	For every multi-index $l= (l_1,\dots l_n)\in\Integer^n_{\geq 0}$, of length $\abs{l} = l_1+\dots l_n$ we write
	$$D_{\sigma, l} \equiv D_{\sigma, 1}^{l_1}\circ\dots\circ D_{\sigma,n}^{l_n}$$
	Then the set $$\left\{v_l(\theta,\sigma) \equiv D_{\sigma, l} (v) (\theta,\sigma) \,\in S(\torus) \, 
	:\, l\in\Integer^n_{\geq 0}\text{,  } v(\theta,\sigma) := e^{-\pi i \krprod{\theta}{\theta}/\sigma}  \right\}$$ is a complete set of eigenvectors for $\hat{\Delta}_{\sigma}$ with corresponding eigenvalues
	$$\hat{\Delta}_{\sigma}v_l(\theta,\sigma) = 2k\left(\abs{l}+ n\frac 12\right)v_l (\theta,\sigma)\text{.} $$
\end{proposition}
\begin{remark}
	This Proposition follows from standard theory of multi-dimensional Hermite polynomials. In fact, the eigenfunctions of $\hat{\Delta}_\sigma$ depend on the the choice of a basis $\{\opb_j\}_{j} $ of $\torus$. If we choose the basis to be \emph{orthonormal} with respect to $\krprod{\cdot}{\cdot}$   the eigenfunctions $\{v_l\}_{l}$ define a Hilbert basis of $\SqInt(\torus)$ with respect to the inner product $\inn{\cdot}{\cdot}_\torus$ from \eqref{eq:inners}.
	Indeed one can easily verify that, in such a basis, $v_l(\theta,\sigma) = (-1)^{\abs{l}}H_{l,\sigma}(\theta) v(\theta,\sigma)$ where 
	\begin{align*}
	&H_{l,\sigma}(\theta) = (-1)^{\abs l}e^{\alpha\krprod{\theta}{\theta}}\partial^l   e^{-\alpha\krprod{\theta}{\theta}},\\
	& \alpha \equiv \pi i \frac{\conj{\sigma}-\sigma}{\sigma\conj{\sigma}}\text{, } \qquad
	 \partial^l = \partial_1^{l_1}\dots\partial_{n}^{l_n}.
	\end{align*}
	which are orthogonal polynomials ({\em multi-dimendsional Hermite polynomials}). Moreover $v(\theta,\sigma)\conj{v(\theta,\sigma) }= e^{-\alpha\krprod{\theta}{\theta}}$, so the $v_l$ are orthogonal $\SqInt(\torus)$ functions. Furthermore recall the following {\em Mehler formula} for every $w\in\Complex$ such that $\re(w)>0$
	\begin{align}
	\nonumber
	\sum_{l\in\Integer_{\geq 0}^n} \frac{w^{\abs l}}{\inn{ v_l}{v_l}_\torus}&H_{l,\sigma}(x)H_{l,\sigma}(y) =\\
	\label{eq:Mehler}
	&\sqrt{\left(\frac{\alpha}{\pi(1-w^2)}\right)^n}
	 \exp\left(\frac{\alpha}{1-w^2}\left(2w\krprod{x}{y}- w^2(\krprod{x}{x}+\krprod{y}{y})\right)\right)
	\end{align}
\end{remark}
Choose $\bpar\in\Complex$ such that $\abs{\bpar} =1$, $\re(\bpar)>0$ and 
\begin{equation}
is = k\frac{1-\bpar^2}{1+\bpar^2}\text{,}
\end{equation}
The choice of a square root of $-\bpar^2$ gives
\begin{align}
&e^{-4kr}= -\bpar^2 &e^{-2kr} =i\bpar
\end{align}
From Proposition \ref{pr:DeltaSpectra} we see that $\hat\Delta_{\sigma} = 2k\left(\hat N_{\sigma}+\frac n2\right)$ where $\hat N_{\sigma}$ has spectrum equal to $\Integer_{\geq 0}$. So, its exponential can be written as 
$$e^{-r\hat\Delta_{\sigma}} = e^{-knr}e^{-2kr\hat N_{\sigma}} =(i\bpar)^{\frac n2} (i\bpar)^{\hat N_{\sigma}}.$$
Finding an explicit expression for $e^{-r\hat\Delta_\sigma} \psi = f$ can be obtained via the kernel $k_{\sigma, \bpar}$ as follows
\begin{align}
&f(\theta,\sigma, \bpar) = \int_{\torus} k_{\sigma, \bpar}(\theta,\tilde{\theta}) \psi(\tilde{\theta})\dvol_k(\tilde{\theta})\text{.}
\end{align}
Recall the Hilbert bases $\{v_l\}_l$ for $\SqInt(\torus)$ diagonalizing $\hat\Delta_{\sigma}$. We can rewrite the kernel as
\begin{align}
k_{\sigma, \bpar} (\theta,\tilde{\theta}) =(i\bpar)^{\frac n2} \sum_{l\in\Integer_{\geq 0}} \frac{(i\bpar)^{\abs l}}{\inn{v_l}{v_l}}_{\torus} v_l(\theta,\sigma)\conj{v_l(\tilde{\theta},\sigma)}
\end{align}
Since $\Re(-\bpar^2) = 1-2(\Re \bpar)^2<0$, Mehler's Formula (\ref{eq:Mehler}) gives the following explicit kernel 
\begin{lemma}\label{lm:Kernel}
	For $i\bpar = e^{-2kr}$ we have
	\begin{align*}
	e^{-r\hat\Delta_{\sigma}}&\psi (\theta) =\\
	&(ib)^{\frac n2}\sqrt{\left(\frac{\alpha}{\pi(1+\bpar^2)} \right)^n}
	\int_{\torus}\exp\left(\frac{\alpha}{1+\bpar^2} \left(2i\bpar\krprod{\theta}{\tilde{\theta}}+ \bpar^2(\krprod{\theta}{\theta}+\krprod{\tilde{\theta}}{\tilde{\theta}})\right) \right)\\
 	&\hspace{5cm} v(\theta,\sigma)\conj{v(\tilde{\theta},\sigma)}\psi(\tilde{\theta}) \dvol_k(\tilde{\theta})
	\end{align*}
	where 
	$$\alpha = \pi i\frac{ \conj\sigma -\sigma}{\sigma\conj\sigma} \text.$$
\end{lemma}

\emph{Proof of Theorem \ref{th:qr}}
The conjugation relation stated in proposition \ref{pr:qr} was already established on equation \eqref{eq:qr}. 
In the case $\sigma = i\bpar$ the operator from Lemma \ref{lm:Kernel} takes the form 
$$ e^{-r\hat{\Delta}_{i\bpar}}\psi(\theta) =  \conj{\bpar}^{\frac 12} e^{\pi(\bpar-\conj{\bpar})\krprod{\theta}{\theta}}\circ\mathcal{F}_\torus(\psi)(\theta)$$
and the explicit expressions follow easily.
$\hfill \square$\\

\section{Remark on Compact Quantum Chern--Simons Theory}\label{sc:Compact}
In this last section we want to make precise contact with the results known for Chern--Simons theory with gauge group $G$.
With reference to \cite{Je}, Proposition $4.2$, consider $\glie$ simply laced.
The matrix elements $S_{\alpha,\beta}$ and $T_{\alpha, \beta}$ presented there for the level-$k$ theory, corresponds to our matrix elements  
$$\hat{\rho}_{k+h,1}(\varphi)_{\frac{(\alpha+\rho)}{k}, \frac{(\beta+\rho)}{k}}\text{ for }\varphi = S\text{ or }T,$$
with the specific choice of $j= i^{-\abs{\Delta_+}}$ and $\omega = e^{\pi i \rprod{\rho}{\rho}_1/h}$. Here $h$ and $\rho$ have the same meaning as the ones in \cite{Je}. The indexes $(\alpha+\rho)/k$ and $(\beta+\rho)/k$  run in the set $I_{(k+h)}$, that corresponds to the weights in the interior of the $(k+h)$-alcove, exactly as in \cite{Je} for the level $k+h$ representation.
The relation between the shift on the level $k\mapsto k+h$ and the corresponding shift on weights $\alpha\mapsto \alpha + \rho$ is explained in \cite{Je}, section A.5.\\
We remark that for $G= \SU(N)$ without the shift $k\mapsto k+h$, we find similar formulae in \cite{Weits}, but there for compact Chern--Simons theory.

\begin{appendices}

\section{Notation from Lie Theory}\label{ap:Lie}
In this appendix we recall some notation and facts from Lie theory that is useful to formulate Chern--Simons theory in genus $1$. We follow the presentations in \cite{ADW,Je} quit closely.\\
Let $G$ be a \emph{compact}, \emph{simple}, \emph{connected}, \emph{simply connected} Lie group, and let $\glie$ be its Lie algebra. Let $\Torus\subset G$ be a maximal torus in $G$ and let $\torus\subseteq \glie$ its Lie algebra, which is a Cartan subalgebra of $\glie$. The \emph{rank} of $G$ is defined as the $\dim\Torus$ and is denoted by $n$.\\
There is a preferred inner product $\rprod{\cdot}{\cdot}_1$ in $\torus$ which is determined, up to a positive scalar $K$, by the trace as $\rprod{X}{Y}_1 = -K\tr(XY)$.\\ 
Associated to $\glie$ we have the \emph{root system} thought of as a certain finite set of $\alpha\in\torus^*$ in the dual of the Cartan subalgebra. Via $\rprod{\cdot}{\cdot}_1$ we can identify $\torus$ and $\torus^*$. In this way the roots can have only two possible lengths.\\
We fix the normalization factor $K>0$ so that $\rprod{\alpha}{\alpha}_1 = 2$ for $\alpha$ a \emph{longest root}. For a \emph{short root} $\beta$ one has $\rprod{\beta}{\beta}_1 = 2/p$, where $p\in\Integer$ depends only on $\glie$. A Lie algebra is called \emph{simply laced} if all its roots have the same length.\\
A fixed set of \emph{positive roots} is denoted $\Delta_+$, while the corresponding set of \emph{simple roots} $\alpha_i$, $i = 1$,$\dots$,$n$, is denoted $\Delta\subset \Delta_+$.\\
To each root $\alpha$ we can associate a \emph{coroot} $h_\alpha\in\torus$ such that $h_\alpha(\alpha) = 2$. Via the preferred inner product we can identify $h_\alpha =  \frac{2\alpha}{\rprod{\alpha}{\alpha}_1}$. In particular the norm of the coroots of a long (resp. short) root $\alpha_l$ (resp. $\alpha_s$) are 
\begin{align}
&\rprod{h_{\alpha_l}}{h_{\alpha_l}}_1 = 2 & \rprod{h_{\alpha_s}}{h_{\alpha_s}}_1 = 2p\text{, for some }p\in\Integer_{>0}
\end{align}
The \emph{coroot lattice} $\coroot\subset \torus$ is the $\Integer$--span of the coroots. This gives the identification
\begin{align}
	\Torus = \torus / \coroot.
\end{align}

For every integer $k\geq 1$ we define a $k$-scaled inner product $\krprod{\cdot}{\cdot}\equiv k\rprod{\cdot}{\cdot}_1$. We use it to define the dual of $\coroot$ with respect to $\krprod{\cdot}{\cdot}$
\begin{align*}
	\kweights \equiv \left\{\gamma\in \torus \text{ such that } \krprod{\gamma}{\lambda} \in\Integer\text{, for every } \lambda\in\coroot \right\}.
\end{align*}
$\coroot$ is integral with respect to $\krprod{\cdot}{\cdot}$, meaning that $\coroot\subset\kweights$ for all $k$. In the special case $k=1$, $\weights = \Lambda_{(1)}^w$ is called the \emph{weights lattice}.\\

We will denote by $\kZcal$ the finite abelian group
\begin{align}
\kZcal\equiv \kweights/\coroot.
\end{align}
For any full rank lattice $\Lambda$ we write $\vol_k{\Lambda}$ for the volume of a fundamental domain of the action of $\Lambda$ on $\torus$ with respect to the inner product $\krprod{\cdot}{\cdot}$. For a finite set $C$, $\abs{C}$ is its cardinality.
We have
\begin{align}
&\Vol_1{\kweights} = \frac{\Vol_1{\weights}}{k^n} = \left(\vol_1{\coroot}k^n\right)^{-1},
&\abs{\kZcal} = \frac{\Vol_1{\coroot}}{\Vol_1{\kweights}} = \abs{\coroot}^2k^n
\end{align}

The Weyl group $W$ is the group generated by the reflections $s_\alpha$ through the hyperplane orthogonal to $\alpha$, i.e.
\begin{align}
s_\alpha(v) = v - 2\frac{\rprod{\alpha}{v}}{\rprod{\alpha}{\alpha}}.
\end{align}
A fundamental domain for the action of $W$ on $\torus$ is called a \emph{Weyl Chamber}.   After a choice of a set of simple roots $\Delta$ is fixed, we can identify a preferred chamber $F_0$, called \emph{fundamental Weyl chamber} as 
\begin{align}
F_0\equiv \left\{ x\in\torus \text{ : } \rprod{x}{\alpha}_1 >0\text{, for all }\alpha\in\Delta \right\}.
\end{align} 
In this way $\torus$ is decomposed into $\abs{W}$ Weyl chambers and $W$ acts simply transitively on the set of them, i.e.
\begin{align}
\torus = \bigcup_{w\in W} \overline{w\left({F_0}\right)}.
\end{align}
The space $\SqInt(\torus,\inn{\cdot}{\cdot}_\torus)$ is defined with respect to the following inner product
\begin{align}
\inn{f}{g}_\torus \equiv \int_{\torus}f(x)\conj{g(x)}\dvol_k(x),
\end{align}
where $\dvol_k$ is the measure on $\torus$ induced by $ \krprod{\cdot}{\cdot}$.\\
An $f\in\SqInt(\torus)$ can be projected orthogonally into the subspace of $W$--\emph{invariant} (resp. $W$--\emph{anti-invariant}) functions $\SqInt_{+}(\torus)$ (resp. $\SqInt_{-}(\torus)$) as follows
\begin{align}
&P_{+}(f) (x) = \frac{1}{\abs{W}}\sum_{w\in W} f(w(x)) &\text{\em (inariant)}\\
&P_{-}(f) (x) = \frac{1}{\abs{W}}\sum_{w\in W} \det(w)f(w(x)) &\text{\em (anti--inariant)}
\end{align}
where $\det(w)$ is the same as the sign character. In particular, both $\SqInt_{+}(\torus)$ and $\SqInt_{-}(\torus)$ are isomorphic to $\SqInt(F_0)$ by restriction to a Weyl chamber.\\
The \emph{affine Weyl group} $W_a$ is the semidirect product $\coroot\rtimes W$. Its \emph{fundamental alcove} $A\subset \torus$ is the set
\begin{align}
A \equiv\left\{ x\in\torus\text{ : } 0<\rprod{x}{\alpha}_1 <1, \text{ for all } \alpha\in\Delta_+ \right\},
\end{align}
and $\overline{A}$ is a fundamental domain for the action of $W_a$ on $\torus$. Denote by $\Omega^R\subset \torus$ the fundamental domain of $\coroot$ containing $A$. We have 
\begin{align}
\vol_1{\coroot} = \vol_1(\Omega^R) = \abs{W}\vol_1(A).
\end{align}

We define the following two sets of indexes 
\begin{align}
& I_k \equiv \kweights\cap A,
&\overline{I}_k \equiv \kweights\cap \overline{A}.
\end{align}
The set $\overline{I}_k$ can be thought as a fundamental domain for the action of $W$ on $\kZcal$ or the action of $W_a$ on $\kweights$. By scaling by $k$  we get a bijection between $\overline I_{k}$ and the weights in the fundamental $k$-alcove (compare with \cite{Je} equation A.16 and section A.5).\\
The space of functions $x\colon \kZcal \longrightarrow \Complex$ is a finite dimensional vector space of dimension $\abs{\kZcal}$. A basis is given by the functions $\delta_\gamma$, $\gamma\in\kZcal$ defined as
\begin{align}
\delta_\gamma({\gamma^0}) = \delta\left(\gamma-\gamma^0 \mod \coroot\right).
\end{align}
The subspaces $\overline{C}_{(k)}$ and $C_{(k)}$ of, respectively, $W$--invariant and $W$--anti-invariant functions are then spanned, respectively, by the following two basis 
\begin{align}
&\overline{C}_{(k)} = \spane_{\Complex} \langle e_\gamma\rangle_{\gamma\in\overline{I}_{k}}
&{C}_{(k)} = \spane_{\Complex} \langle \tilde e_\gamma\rangle_{\gamma\in{I}_{k}}
\end{align}
\begin{align}
\label{eq:InvFinBasis}
&e_{\gamma} \equiv \frac{1}{\sqrt{\abs{W}}}\sum_{w\in W} \delta_{w(\gamma)}, 
&\text{for every }\gamma\in \overline{I}_k,\\
\label{eq:AntiFinBasis}
&\tilde e_{\gamma} \equiv \frac{1}{\sqrt{\abs{W}}}\sum_{w\in W}\det(w) \delta_{w(\gamma)},
&\text{for every }\gamma\in {I}_k.
\end{align}
We use the sets $\overline{I}_{k}$ and $I_k$ to index the  basises $\{e_\gamma\}$ and $\{\tilde{e}_\gamma\}$ respectively.\\
These basises are orthonormal with respect to the product
\begin{align}
\inn{x}{y}_{\kZcal} \equiv \sum_{\gamma\in\Zcal}x(\gamma)\conj{y(\gamma)}.
\end{align} 
\end{appendices}

\bibliographystyle{plain}

\end{document}